\documentclass[oneside,draft,letterpaper,12pt]{amsart}
\usepackage{amsthm}
\usepackage{amsmath}
\usepackage{amssymb}
\theoremstyle{definition}
\def\bd#1#2#3{{#1}_{#3}({#2})}
\def\bn#1#2{B_{#2}({o_{#1}})}
\def\R{\text{$\mathbb{R}$}}

\def\fis#1{\dot{#1}}
\def\lra{\longrightarrow}
\def\pt#1#2{\tau^{#2}_{#1}}
\def\fis#1{\dot{#1}}
\def\lra{\longrightarrow} 
\def\nb{\nabla}

\def\d1#1#2{\frac{d#1}{d#2}}
\def\nbp#1#2{\nabla_{\frac{\partial}{\partial #1}}#2}
\def\tpm{T_p M}
\def\pt#1#2{\tau^{#2}_{#1}}
\def\e#1{\exp_{#1}}
\def\bd#1#2#3{{#1}_{#3}({#2})}
\def\bn#1#2{B_{#2}({0_{#1}})}
\def\p1#1#2{\frac{\partial #1}{\partial #2}}
\def\sfv#1#2#3{\nabla_{\frac{\partial}{\partial#1}}\nabla_{\frac{\partial}{\partial #2}}#3}
\def\to{t_o}
\def\pr{\parallel}
\def\cal#1{\mathcal{#1}}
\def\fD{f_{\Delta}}
\def\fd{f_{\delta}}

\def\ptn#1#2{\chi^{#2}_{#1}}

\def\D{D^2 E(c)}
\def\so{\rightarrow}

\def\part{a=\to \leq t_1 \leq \ldots \leq t_{n-1} \leq t_{n}=b}

\def\D{ D^2 E (\gamma) }
\def\ov{\overline}
\newcommand{\dem}{{\bf Proof:}}
\newcommand{\ed}{\end{document}}
\newcommand{\n}{\noindent}
\newcommand{\We}{A_{\xi}}
\def\Ex{\mathrm{Exp}^{\perp}}
\def\taT{\tilde{T}}

\def\N{\text{$\mathbb{N}$}}
\def\H{\text{$\mathbb{H}$}}

\def\vai{\rightarrow}

\def\mh{M^{H}}
\def\mhe{M^{H-\epsilon}}
\def\tr{(\gamma_1,\gamma_2,\gamma_3)}
\def\tro{(\overline{\gamma_1},\overline{\gamma_2},\overline{\gamma_3})}

\def\gga{(\gamma_1,\gamma_2,\alpha)}
\def\ggao{(\overline{\gamma_1},\overline{\gamma_2},\overline{\alpha})}

\def\trh{(\gamma_1,\gamma_2, \frac{\pi}{2})}
\def\trho{(\overline{\gamma}_1,\overline{\gamma}_2, \frac{\pi}{2})}
\def\ga#1{\gamma_{#1}}
\def\gao#1{\overline{\gamma}_{#1}}
\def\alf#1{\alpha_{#1}}
\def\alfo#1{\overline{\alpha}_{#1}}
\def\lo#1{\overline{l}_{#1}}

\def\si#1{\sigma_{#1}}

\def\ti#1{\tau_{#1}}

\def\bi#1{\beta_{#1}}

\def\pr{\parallel}
\title[Infinite dimensional geometry]{Some results on infinite dimensional 
Riemannian geometry}
\author[L. Biliotti]{Leonardo Biliotti}
\thanks{Research partially supported by CNPq  (Brazil) }
\begin{document}
\theoremstyle{plain}\newtheorem{thm}{Theorem}
\theoremstyle{plain}\newtheorem{prop}[thm]{Proposition}
\theoremstyle{plain}\newtheorem{lemma}[thm]{Lemma}
\theoremstyle{plain}\newtheorem{cor}[thm]{Corollary}
\theoremstyle{definition}\newtheorem{defini}[thm]{Definition}
\theoremstyle{remark}\newtheorem{remark}[thm]{Remark}
\theoremstyle{plain} \newtheorem{assum}[thm]{Assumption}
\theoremstyle{definition}\newtheorem{ex}{Example}
\begin{abstract}
In this paper we will investigate the global properties of complete Hilbert
manifolds with upper and lower bounded  sectional curvature. We shall  prove
the Focal Index lemma that will allow us to extend some classical results of 
finite dimensional Riemannian geometry as Rauch and Berger 
theorems and the Topogonov 
theorem in the class of manifolds in which the Hopf-Rinow theorem holds.
\end{abstract}
\maketitle
\section{Introduction}
\mbox{}
In infinite dimensional geometry, the most of the local results follow from 
general arguments analogous to those in the finite dimensional case 
(see \cite{Kl} or \cite{La}). The investigate of global properties is quite 
hard as finite dimensional case and the theorem of Hopf-Rinow is generic 
satisfy on complete Hilbert manifolds (see \cite{Ek}).
Moreover, the exponential map may not be surjective also when the manifold
is a complete Hilbert manifold. In section 1 we briefly discuss the
relationship between completeness and geodesically completeness (at some
point) and we note  that this facts are equivalent either when a manifold has 
constant sectional curvature or no positive sectional curvature. We conclude 
this section to prove that a group of bijective isometry coincide with
the set of the maps that preserve the distance.   

The fundamental tools
to studying the geometry and topology of the finite dimensional manifolds are 
the Rauch and Berger theorems. These theorems allow us to 
understand the distribution of conjugate and focal points along
geodesic and the geometry of the complete Hilbert manifolds with bounded
sectional curvature. We recall briefly the notion of focal point:
let $N$ be a submanifold of a riemannian manifold $M.$ The exponential map
of $M$ is defined on an open subset $W \subset TM$ and we restrict its on
$W \cap T^{\perp} N,$ that we will denote by $\Ex.$  
A focal point is a singularity of  $\Ex :W \cap T^{\perp} N \lra M  .$
In the infinite dimensional manifolds two species of singularity can be appear:
when the differential  of $\Ex$ fails to be injective (monofocal)
or when the differential of $\Ex$ fails to be surjective (epifocal). Clearly,
when $N=p$ we have exactly the notion of conjugate points. 
In section $2$, we shall study the singularity of the exponential map of 
differential point of view, as in \cite{Ka}, and we shall prove that
always monofocal point 
implies epifocal but not conversely and the distribution of
epifocal points and monofocal points can have cluster points. 
Moreover we will deduce a weak form  of the Rauch theorem.

In section $4$ we prove the fundamental tool to prove the main results:
the {\em Focal Index lemma}. Firstly, we will prove the same version
of the finite dimensional theory; then we note that we shall prove its
in the case when we have only a finite number of epifocal points which
are not monofocal (pathological points). 
Using the above results we will prove the Rauch and Berger
theorems in infinite dimensional geometry, 
when we have at most a finite number of pathological points along a finite
geodesic. The main applications will appear in the last
section and the main result is the Topogonov Theorem in the 
class of the complete Hilbert manifolds on which the Hopf-Rinow theorem holds. 
The prove is almost the same as in \cite{CE} because Rauch, Berger and 
Hopf-Rinow theorems hold. Moreover, 
we shall prove the Maximal Diameter Theorem.
and   two version of sphere theorems, 
with the strong assumption on injectivity 
radius, with pinching $\sim \frac{3}{4}$ and $\frac{4}{9},$
in the class of Hopf-Rinow manifolds. Other simple applications are
two results like Berger-Topogonov theorem, one using Rauch Theorem
and one using Topogonov Theorem, and a results about the image of
the exponential map of a complete manifold with upper bounded sectional 
curvature. Some basic references for
infinite dimensional geometry are  \cite{La} and
\cite{Kl}.
\section{Preliminaries}
\mbox{}
In this section we give some general results of infinite dimensional 
Riemannian geometry and we briefly discuss some of the differences from the 
finite dimensional case. 
We begin recalling some basic  facts and establishing our notation. 

Let $(M,g)$ be an Hilbert manifold modeled on a infinite dimensional Hilbert 
space $\H.$ Throughout this paper we shall assume that $M$ is a connected,
paracompact and Hausdorff space.   
Any tangent space $\tpm$ has a scalar product $g(p),$
depending smoothly on $p,$ and defining on $\tpm$ a norm equivalent to the 
original norm on $\H.$ Using $g,$ we can define the length of piecewise 
differential curve
and it easy to check that for any two points of $M$ there exists
a piecewise  differential curve joining them . Hence, we can introduce a 
metric $d,$ defining $d(x,y)$
the infimum of the lengths of all differential paths joining $x,y$ and one 
can prove that  $(M,d)$ is a metric space (see \cite{Pa} ) and $d$ 
induces the same topology of $M.$ 
As in the finite dimensional case, $M$ admits a unique 
Levi Civita connection $\nb$ (see \cite{La}) defined by the 
Koszul's formula. We recall that the criterion of tensoriality in infinite
dimensional geometry doesn't hold (see \cite{BMT}) so we must deduce all 
properties of $\nb$ by its local expression. 

Let $c:[a,b] \lra M$ be 
a smooth curve. For any $ v \in T_{c(a)} M$ there exists a unique vector 
field $V(t)$ along  $c$ such that $ \nb_{\fis c (t)}{V(t)}=0.$ 
Moreover, the Levi Civita connection satisfies
$
\nb_X Y (p) = \frac{d}{dt} \mid_{t=0}  \pt{t}{0} (Y),
$
where $\pt{t}{0}$ is a parallel transport along any curve $\gamma$ such that
$\gamma (0)=p$ and $\fis{\gamma} (0)= X(p).$ 

A geodesic in $M$ is a smooth curve $\gamma$ which satisfies 
$\nb_{\fis{\gamma}} \fis{\gamma}=0.$ Using the theorems of existence, 
uniqueness and smooth dependence on the initial data, we may prove the 
existence, at any point $p,$ of the exponential map $\e p,$ that is defined 
in a neighborhood of the origin in $\tpm$ by setting $\e p (v)=\gamma (1),$
where $\gamma$ is the geodesic in $M$ such that $\gamma (0)=p$ and 
$\fis{\gamma }(0)=v.$ This map is smooth and is a local diffeomorphism,
$d( \e p )_{0}=id,$ in a neighborhood of the origin in $\tpm$ by the inverse 
function theorem. Moreover, there exists an open neighborhood $W$ of  
$\{ 0_p \in \tpm: \ p \in M  \}$ in the tangent bundle $TM,$ 
such that the application 
$\exp (X_p)= \e p (X)$ is defined and differentiable. 

Generally, the  local theory works as in 
finite dimensional geometry and results as the 
Gauss lemma and existence of convex neighborhoods hold in infinite 
dimensional Riemannian geometry and the curvature tensor is defined as
follows: let 
$x,y,z \in \tpm,$ we extend them to vector fields $X,Y,Z$ and define
$R(x,y)z=\nb_X \nb_Y Z - \nb_Y \nb_X Z - \nb_{[X,Y]} Z.$
It is easy to check that $R$ doesn't depend on the extension and it is 
antisymmetric in $x,y$ and satisfies the first Bianchi identity
$
R(x,y)z + R(z,x)y + R(y,z)x =0.
$
Given any plane $\sigma$ in $\tpm$ and let $v,w \in \sigma$ be 
linearly independent. We define the sectional curvature 
$K(\sigma)$ to be
\begin{small}
$$
\frac{g(R(v,w)w,v)}{g(v,v) g(w,w) - g(v,w)^2}.
$$
\end{small}
It is easy to check that $K$ doesn't depends on the choice of the 
spanning vectors and the curvature tensor $R$ is 
completely determined by the sectional curvature.   
Moreover, as in the finite dimensional case, we shall prove the 
Cartan
theorem, see \cite{Kl} page 114, 
in which the existence of a local isometry is characterized by a 
certain property of the tensor curvature. 

The global theory on  Hilbert manifolds is more difficult: for 
example the Hopf-Rinow theorem fails. 
We recall that a Hilbert manifold is called complete if $(M,d)$ 
is complete as a metric space. On the other hand, we say that a manifold
$M$ is geodesically complete  at a point $p$ if the 
exponential map is defined
in $\tpm$ and $M$ is geodesically complete if it is geodesically complete
for all point $q \in M.$ It is easy to check the following 
implication: \\
 
\n
$M$ complete $\Rightarrow$ $M$ is geodesically complete   
$\Rightarrow$ $M$ is geodesically complete at some point. \\

If $M$ is a finite dimensional manifold, the above sentences are
equivalent, thanks to the Hopf-Rinow theorem. 
Grossman, see \cite{Gr} constructs a simply connected complete Hilbert 
manifold on which there exist two points  which cannot be  
connected with a minimal geodesic but the exponential map is surjective.
On the other hand Ekeland (\cite{Ek})
proved that the Hopf-Rinow theorem is generically satisfied, i.e. if one takes
a point $p \in M,$ in a complete Hilbert manifold $M,$ 
the set of  points $q\in M$
such that there exists a unique minimal geodesic joining $p$ and $q$ is a
$G_{\delta}$ set and in particular it is a dense subset of $M.$
The Hopf-Rinow theorem implies also that the exponential map must be surjective
on complete finite dimensional Riemannian manifolds. Atkin (see \cite{At}) 
showed that there exists a complete Hilbert manifold $M$ on which,  
at some point $p \in M,$ 
$\e p (\tpm)$ is not surjective. 
Let $q \in M - (\e p (\tpm)).$ Then $M-\{ q \}$ is not 
complete as metric space, but, clearly, is geodesically complete at $p.$
In particular: 
{\em in infinite dimensional Riemannian geometry, geodesically complete at some
point doesn't imply completeness.} Moreover, Atkin (see \cite{At1}) 
constructed some infinite dimensional Hilbert manifolds in which the induced
metric is incomplete, but geodesically complete and any two point may be 
joined by minimizing geodesic.  
The above discussion justify the following definition.
\begin{defini}
A complete Hilbert manifold $M$ is called Hopf-Rinow if for every $p,q \in M$ 
there exists at least a minimal geodesic joining $p$ and $q.$
\end{defini}
In \cite{El} El\'{\i}ason showed that the Sobolev manifolds, i.e. 
the spaces of the Sobolev 
sections of a vector bundle on a compact manifold, are Hopf-Rinow. Other
class of manifolds which are  Hopf-Rinow are  the simply connected 
complete Hilbert manifolds such that their sectional 
curvature is non positive;
indeed  the Cartan-Hadamard theorem holds, see \cite{Gr} and \cite{Mca}. 
Furthermore, see \cite{La}, 
we may prove that an Hilbert manifold $M$ with non positive sectional 
curvature is a
complete Hilbert manifold if and only if $M$ is geodesically complete at
some point. It is easy to check that the same holds for a
manifolds with constant sectional curvature: indeed, using the same arguments
used in 
the classification of the simply connected complete Hilbert manifolds with 
constant sectional curvature, we may prove our claim.

The Bonnet theorem was proved by Anderson, see \cite{An}; however we cannot 
conclude any information about fundamental group since we may prove,
see \cite{B}, that there exist 
infinite groups acting isometrically and properly
discontinuously on the infinite unit sphere. In particular, Weinstein 
theorem fails: the following example gives an isometry of the unit
sphere of a separable Hilbert space $S(l_2)$ without fixed points, such that
$inf \{ d(x,f(x)) : x \in S(l_2) \} =0.$
$$
\begin{array}{lcl}
f(\sum_{i=1}^{\infty} x_i e_i)&=&
\sum_{i=1}^{\infty} (\cos(\frac{1}{i}) x_{2i-1} \ + \
                     \sin(\frac{1}{i}) x_{2i}) e_{2i-1} \\
                    &+&
                    \sum_{i=1}^{\infty}
                    (\cos(\frac{1}{i}) x_{2i})\ - \
                    \sin(\frac{1}{i}) x_{2i-1})e_{2i}. \\
\end{array}
$$  
We conclude this section  proving that the group of bijective
isometry coincide with the set of 
applications that  preserve the distance.
\begin{prop} \label{isometrie}
Let $F:(M,g) \lra (N,h)$ be a surjective map. Then $F$ is an isometry if and 
only if $F$ preserves the distance, i.e. if and only if $d(F(x),F(y))=d(x,y).$
\end{prop}
The same proof of the finite dimensional case holds, (\cite{KN})  
once we prove the 
following result
$$\lim_{s \so 0} \frac{d(\e p (sX) , \e p  (sY)}{\pr sX - sY \pr}=1.
$$ 
\dem $\ $ let  $r>0$ such that
$\e p: \bn{p}{r} \rightarrow \bd {\cal B}{p}{r},$
between the balls with their respectively metrics, 
is a onto diffeomorphism. Now recall that $d( \e p)_0 =id,$ hence
there exists $\epsilon >0$ and  $0< \eta \leq r$ such that
\[
1 -  \epsilon \leq \pr d (\e p)_q \pr \leq 1  +  \epsilon,
\]
for every $q \in \bn{p}{\eta}$. Let
$m,a \in \bd{\cal B}{p}{\frac{\eta}{4}}$.  By assumption, we 
shall calculate the distance
from $m$ to $a,$ restricting ourself to the curve $c$ on
$\bd{\cal B}{p}{\eta}$.  Then $c(t)=\e p(\xi(t))$ and
\[
(1 - \epsilon) \int_0^1 \parallel \fis {\xi}(t)  \parallel_p dt
\leq \int_0^1 \parallel \fis c (t) \parallel dt \leq 
(1 + \epsilon) \int_0^1 \parallel \fis {\xi}(t)  \parallel_p dt.
\]
That is: for every $\epsilon>0$ there exists
$s_o=\frac{\eta}{4}$ such that, for every $s<s_o$ we have
\begin{small}
\[
 (1 - \epsilon)(\pr sX  -  sY \pr_p ) \leq d(\e p (sX), \e p (sY)) \leq 
 (1 + \epsilon)( \pr sX  -  sY \pr),
\]
\end{small}
that implies our result. QED
\section{Jacobi Flow}  
\mbox{}
The linearized version of the geodesic equation is the 
famous Jacobi equation. In this section, we shall study 
the Jacobi field from the dif\-fe\-ren\-tial point of view getting 
some informations of the 
dis\-tri\-bu\-tion of singular points of exponential map. 
Throughout this section, 
all es\-ti\-mates are formulated in terms of unit speed geodesics 
because one can easily reparametrize: if $J$ is a Jacobi fields 
along $c,$ then 
\begin{center}
$J_r(t)= J(rt)$  is  a Jacobi field along $ c_r (t) = c(rt)$; \\
with $J_r (0)=J(0)$ and ${J'}_r(0)=r J' (0), {c'}_r (0)=r c' (0)$ 
\end{center}   
\begin{defini}
Let $c:[0.a] \lra M$ be a geodesic. A vector field along $c$ is called Jacobi
field if it satisfies the Jacobi differential equation
$$
\nbp{t} \nbp{t}{J}(t) + R(J(t), \fis c (t) ) \fis c (t) =0,
$$
where $\nbp{t}{}$ denotes the covariant derivation along $c.$
\end{defini}
The Jacobi equation is a second order differential equation and, by 
theorems of differential equation in Banach space, see \cite{La},
the solutions are defined in the 
whole domain of definition and the set of Jacobi
fields along $c$ is a vector space  isomorphic to 
$T_{c(0)} M \times T_{c(0)} M$ under the map 
$J \lra (J(0), \nbp{t}{J} (0) ).$   We recall also that the Jacobi fields
are characterized as in\-fini\-te\-si\-mal variation of $c$ by geodesic. 
We first give a lemma due to Ambrose (see \cite{La} page 243).
\begin{lemma}
Let $c:[0,a] \lra M$ be a geodesic and let $J$ and $Y$ be two Jacobi vector 
fields along $c.$ Then 
$$
\langle \nbp{t}{J} (t), Y(t) \rangle - \langle \nbp{t}{Y} (t), J(t) \rangle
={\rm constant} C.
$$
\end{lemma}
Now, let $c: [0,b] \lra M$ be a geodesic. We will denote by $p=c(0)$ and 
by $\pt{s}{t}$ the 
parallel transport along $c$ between the points $c(s)$ and $c(t)$. We define   
$$
R_s:T_{c(0)} M   \lra T_{c(0)} M, \ \
R_s(X)=\pt{s}{0}(R(\pt{0}{s}(X),\fis{c}(s))\fis{c}(s))
$$
that is a family of symmetric operators in $\tpm$. Take 
$\H_o$ a closed subspace of $\tpm$  and let 
$
A:\H_o \lra \H_o
$
be a continuous and symmetric linear operator. Clearly
$\tpm=\H_o \oplus \H_o^{\perp}$ and we will study the solutions of the
following linear differential equation
\[
\left \{ \begin{array}{l}
T''(s)\ +\ R_s(T(s))=0;\\
T(0)(v,w)=(v,0),\ T'(0)(v,w)=(-A (v),w),
\end{array}
\right.
\]
that we call {\em Jacobi flow } of $c.$
Firstly, note  that the family of bilinear applications
$$
\begin{array}{lcl}
& \tpm \times \tpm \stackrel{\Phi(t)}{\lra} \R & \\
& (u,v) \lra  \langle T(t)(u),T'(t)(v) \rangle & \\
\end{array}
$$
is symmetric; indeed is symmetric in $t=0,$ because $A$ is a symmetric 
operator, and 
\[
( \langle T(t)(u),T'(t)(w) \rangle \ - \  \langle T(t)(w),T'(t)(u) \rangle )'
\]
is zero. The solutions of the Jacobi flow are exactly the Jacobi fields.
\begin{prop}
Let $(v,w) \in \tpm=\H_o \oplus \H_o^{\perp}.$ 
Then the  Jacobi field with 
$J(0)=v \in \H_o$ and $\nbp{t}{J(0)}\ + \  A (J(0))=w \in {\H_o}^{\perp},$
is given by  $\pt{0}{t}(T(t)(v,w))$.
\end{prop}
\n
\dem $\ $
Let  $Z(t)$ be a parallel transport of $u\in \tpm$ along $c.$  
We indicate by $Y(t)=\pt{0}{t}(T(t)(v,w))$. Hence
\begin{eqnarray}
\langle \sfv{t}{t}{Y(t)},Z(t) \rangle  &=& \langle Y(t),Z(t) \rangle''   
\nonumber \\
                       &=& \langle T(t)(v,w),u \rangle ''    \nonumber \\
                       &=& - \langle R(Y(t),\fis {c}(t)) \fis{c} (t),Z(t) 
\rangle. \ \mathrm{QED}
\nonumber
\end{eqnarray}
Our aim is to study the distribution of singular points of the 
exponential map along a geodesic
$c:$ then it is very useful compute the adjoint operator of $T(b).$ 
Let  $u \in \tpm$ and let $J$ be the Jacobi field along the geodesic 
$c$ such that
$
J(b)=0, \ \nbp{t}{J(b)}=\pt{0}{b}(u).
$
By the Ambrose lemma we have 
\[
(1)\ \langle T(b)(v,w),u \rangle = \langle T(0)(v,w),
\nbp{t}{J}(0) \rangle \ - \ \langle
T'(0)(v,w),J(0) \rangle.
\]
We denote by $\ov{c}(t)=c(b-t)$ and let 
\[
\left \{ \begin{array}{l}
\taT''(s)\ +\ R_s(\taT(s))=0;\\
\taT(0)=0,\ \taT'(0)=id,
\end{array}
\right.
\]
be the Jacobi flow of $\ov{c}:[0,b] \lra M$. 
It is easy to check that  $J$ is a vector field along $c,$ then 
$\overline{J}(t)=J(b-t)$ is the Jacobi filed along 
$\ov{c}$ such that
$\nbp{t}{\overline{J}}(b)=-\nbp{t}{J}(0)$. Then $(1)$ becomes
$$
\begin{array}{lcl}
\langle T(b)(v,w),u \rangle &=& \langle (v,o),
\pt{b}{0}(\taT'(b)(-\pt{0}{b}(u))) \rangle \\
&-&  \langle (-A(v),w),\pt{b}{0}(\taT (b)(-\pt{0}{b}(u)) \rangle \\
\end{array}
$$
and the adjoint operator is given by
\begin{small}
$$\begin{array}{l}
 \langle T^*(b)(u),(v,0) \rangle  = -  \langle \pt{b}{0} 
 (\taT'(b) (\pt{0}{b}(u) ) ) ) +
  A( p_t ( \pt{b}{0}(\taT(b) ( \pt{0}{b}(u) ) ) ) ) ,(v,0) \rangle ,   \\
\langle T^*(b)(u), (0,w) \rangle =\
\langle \pt{b}{0}(\taT(b)(\pt{0}{b}(u))), (0,w) \rangle .    \\
\end{array}$$
\end{small}
where $p_t$ is the component along $\H_o.$
\begin{prop} \label{p1} 
There exists a bijective correspondence between the kernel of $ T(b)$ and 
the kernel of $T^* (b).$ 
\end{prop}
\n
\dem 
$\ $ let $w \in \tpm$ such that $T(b)(w)=0.$ Then the Jacobi field 
$$
Y(t)=\pt{0}{t}(T(t)(w))=\pt{b}{t}(\taT(t)(\overline{w}))
$$
and  $\overline{w} \in T_{c(b)}M$ is unique. 
Using the boundary condition
on $Y(t)$ we have $T^* (b) (\pt{b}{0} (\ov{w}))=0.$ Vice-versa, if 
$T^* (b) (w)=0$ it is easy to check that 
$\pt{b}{t} \taT (t) (\pt{0}{b}(w))= \pt{0}{t} (T(t) (\ov{w})$ for some
$\ov{w} \in \tpm.$ Clearly
$T(b) (\ov{w})=0$ and $w$ can be obtained, starting from $\ov w$, with
the above arguments. QED\\

We recall that
$Ker T^* (b)=\ov{Im T}^{\perp},$ hence we have proved also 
the following result
\begin{prop} \label{singular1}
If $T(b)$ is injective then  $ImT(b)$ is a dense subspace.
\end{prop}
We will study the behaviour of the Jacobi flow either when 
$\H_o=0$ or $\H_o= \tpm$ and $A=0.$ We will denote by 
$f_{\Delta}$ a solution of the differential equation
$f''(t)\ +\ \Delta(s)f(s)=0$ with
$f(0)=0$, $f'(0)=1$, when $\H_o=0,$ or
with $f(0)=1$, $f'(0)=0$, when $\H_o=\tpm$ and $A=0$.
Firstly, we note that in \cite{Ka} it is proved the following results that
holds in our context. 
\begin{prop} \label{p2}
Let  $\delta(s) \leq
\langle R(u,\fis c(s))\fis c(s),u \rangle \leq \Delta(s)$ be an upper and
lower curvature along $c$. Hence
\begin{enumerate}
\item if $T(s)$ is invertible then 
\[
( \langle T'T^{-1} u,u \rangle )' \leq -(\Delta(s)\ + \
\langle T'T^{-1}u,T'T^{-1}u \rangle^2);
\]
\item $\pr T(s)u \pr \geq \fD(s) \langle u,u \rangle^{\frac{1}{2}},$
$\fD(s)$ being positive in  $0< s \leq s_o$;
\item $\langle T'u,Tu \rangle \fD \geq  \langle Tu,Tu \rangle  \fD'$, 
$\fD (s)$ being positive in $0 < s \leq s_o$;
\item $ \langle T'u,Tu \rangle \fd \leq  \langle Tu,Tu \rangle \fd'$,
$0< s \leq s_o$ if $T(s)$ is invertible in   $0<s \leq s_1$ $(s_1 \geq s_o)$;
\item if $T(s)$ is invertible in  $0< s \leq s_1$ then
\[
\pr T(s) \pr \fd(t) \geq \pr T(t) \pr \fd (s), \ 0 \leq s \leq t
\leq s_1;
\]
\item $\pr T(s) \pr \leq \fd(s) <u,u>^{\frac{1}{2}},\ 0 \leq s \leq s_1$.
\end{enumerate}
\end{prop}
\begin{cor} \label{singular2}
Let
\[
\left \{ \begin{array}{l}
T''(s)\ +\ R_s(T(s))=0; \\
T(0)(v,w)=(v,o),\ T'(0)(v,w)=(-A (v),w),
\end{array}
\right.
\]
be the  Jacobi flow. We suppose that we have a upper bound sectional 
curvature, i.e. $K \leq H. $ Then:
\begin{itemize}
\item [{\bf (a)}]
if $\H_o=0$ then the Jacobi flow is a topological isomorphism for
$t \geq 0$, if $H \leq 0,$ 
or for  $ 0 \leq t < \frac{\pi}{\sqrt{H}},$ if $H>0$;
\item [{\bf (b)}]
if $\H_o=\tpm$ and $A=0$ then then Jacobi flow is a topological isomorphism
for $t \geq 0$, if $H \leq 0,$ or
for $ 0 \leq t < \frac{\pi}{2 \sqrt{H}},$ if $H>0$.
\end{itemize}
\end{cor}
\n
\dem $\ $ using property $(2)$ of proposition \ref{p2} 
and proposition \ref{p1} we have our claim.\\

Another interesting corollary is a weak form of the Rauch theorem.
\begin{thm} \label{rauch} {\bf (Rauch weak)} 
Let  $c:[0,b] \lra M$ be a unit speed geodesic. Suppose that we have a lower 
and upper bound of the sectional curvature, i.e.  
$L \leq  K(\fis c(t), v) \leq H,$
for any $t$ and  $v \in T_{c(t)}M$ such that 
$\langle v,\fis c (t) \rangle=0,\ \langle v, v \rangle =1.$
Then $d (\e{p})_{t \fis c (0)}$ is a topological isomorphism for every
$t \geq 0$, if $H \leq 0$ and for every
$0 \leq t < \frac{\pi}{\sqrt{H}}$ if $H>0$. Moreover, 
in the above neighborhoods, we have
\[
 \frac{f_{H}(t)}{t} \pr v \pr  \leq
 \pr d ( \e p )_{t \fis c (0)} (v) \pr \leq
 \frac{f_{L}(t)}{t} \pr v \pr .
 \]
\end{thm}
\section{Focal Index lemma and the Rauch and Berger theorems}
\mbox{}
Let $N$ be a submanifold of an Hilbert manifold $(M,g),$ i.e. $N$ is an 
Hilbert manifold and the inclusion $i: N \hookrightarrow M$ is an embedding.  
Let $\gamma:[a,b] \lra M$ be a geodesic such that
\begin{enumerate}
\item $\gamma(a)=p \in N$;
\item $\fis{\gamma}(a)=\xi  \in {T_p N}^{\perp}$.
\end{enumerate}
\begin{defini}
A Jacobi vector field along $\gamma$ is called  $N$-Jacobi if it satisfies
the following boundary conditions:
\[
Y(a) \in T_p N,\ \nbp{t}{ Y(a)} \ +\ A_{\xi}(Y(a)) \in T_p^{\perp} N,
\]
where $\We$ is the operator of  Weingarter relative to $N$.
\end{defini}
The Jacobi flow along $\gamma$ on which $A= \We$ is called Jacobi flow
along $\gamma$ of $N.$ Let $W$ be the open subset of $TM$ on which $\exp$ is
defined.
We denote by  $\Ex :T^{\perp} N \cap W \lra M,$ 
$\Ex (X)=\exp^{M}(X).$ It is
easy to prove  that 
$  \mathrm Ker ( d ( \Ex )_{\to \xi} )$ is isomorphic to 
the
subspace of the $N$-Jacobi field along the geodesic 
$\gamma (t) = \Ex (t \xi),$ $ t \in [0,t_o],$  
which is zero in $\gamma (t_o).$ 
On the other hand, in
infinite dimensional manifolds, two species of singular points can appear;
so the following definition is justify. 
\begin{defini}
A element $q=\gamma(\to)$ along  $\gamma$ is called:
\begin{enumerate}
\item monofocal if $\mathrm d ( \Ex )_{\to \xi} $
fails to be injective;
\item epifocal if $\mathrm d ( \Ex )_{\to \xi} $
fails to be surjective.
\end{enumerate}
\end{defini}
By proposition \ref{p1} and by the formula of the adjoint of the Jacobi
flow,  we have the following results.
\begin{prop} \label{focal}
Let $N$ be a submanifold of $M$ and let $\gamma :[0,b] \lra M$ such that
$\gamma (0)=p \in N$ and $\xi= \fis{ \gamma} (0) \in {T_p N}^{\perp}.$ Then
\begin{enumerate}
\item if $\gamma (t_o)$ isn't monofocal then the image of
 $\mathrm d( \Ex )_{t_o \xi}$ is a dense subspace;
\item $\gamma (t_o)$ is monofocal then $\gamma (t_o)$ is epifocal;
\item $q=\gamma (t_o)$ is monofocal along $\gamma$ if and only if $p$ is 
monofocal along $\ov{\gamma}(t)=\gamma(t_o-t);$
\item $q=\gamma (t_o)$ is epiconjugate and the image of  
$\mathrm d ( \Ex )_{\to \xi} $ is closed then $p$ is monofocal along
$\ov{\gamma}(t)=\gamma(t_o-t);$ 
\end{enumerate}
\end{prop}
In the degenerate case, i.e.  $N=p,$  then we call  $q=\gamma(\to)$ either
{\em monoconjugate} or {\em epiconjugate} along $\gamma$.
If there exist neither monofocal (monoconjugate)
nor epifocal (epiconjugate)  points then we will say that 
there aren't  focal (conjugate) points along $\gamma$. 
The distribution of singular points  of the exponential map along a finite 
geodesic is different
from the finite dimensional case. Indeed, Grossman showed how the distribution 
of monoconjugate points be able to have cluster points. 
The same pathology appears in the case of focal points  and we shall give a 
pathological example of the   distribution of monofocal and epifocal 
points along a  finite geodesic. 
\begin{ex} \label{esempio}
{\rm
Let
$M= \{ x \in l_2 :\ x^2_1\ +\ x^2_2 \ + \ \sum_{i=3}^{\infty} a_i x^2_i=1 \}$,
where $( a_i )_{i \in \N}$ is a positive sequence of real number. 
It is easy to check that
\[
\gamma(s)= \sin (s) e_1 \ + \ \cos (s) e_2
\]
is a geodesic and   $T_{\gamma(s)}M=<\fis{\gamma}(s),e_3,e_4,\ldots >$. 
Let $N$ be a geodesic submanifold defined by $\fis{\gamma} (0).$ 
We shall restrict ourself to the normal Jacobi fields.
We note that for $k \geq 3$
\[
E_k:= \{ x^2_1\ +\ x^2_2 \ +\ a_k x^2_k =1\ \} \hookrightarrow M
\]
is totally geodesic; then $K(\fis{\gamma}(s),e_k)=a_k$ and
the Jacobi fields, with boundary conditions
$J_k(0)=e_k$, $\nbp{t}{J_k(0)}=0$, are given by 
$J_k(t)=\cos (\sqrt{a_k}t)e_k$. Hence
\[
d (\Ex)_{s \fis{\gamma}(0)} (\sum_{k=3}^{\infty} b_k e_k )=
\sum_{k=3}^{\infty} b_k \cos (\sqrt{a_k}s) e_k.
\]
Clearly, the points $\gamma (r_k^m)$, $r_k^m=\frac{m \pi}{2 \sqrt{a_k}}$ 
are monofocal. Specifically, let 
 $a_k=(1- \frac{1}{k})^2.$ The points $\gamma (s_k)$,
$s_k=\frac{k \pi}{2 (1-k)}$ are monofocal, $s_k \so \frac{\pi}{2}$ and
\[
d (\Ex)_{\fis{\gamma}(\frac{\pi}{2})) } ( \sum_{k=3}^{\infty} b_k e_k)=
      \sum_{k=3}^{\infty} b_k  \cos (\frac{k-1}{k} \frac{\pi}{2}) e_k.
\]
In particular, $\gamma(\frac{\pi}{2})$ is not monofocal along $\gamma.$ 
On the other hand if $\sum_{k=3}^{\infty} \frac{1}{k} e_k=
d (\Ex)_{\gamma(\frac{\pi}{2})}(\sum_{k=3}^{\infty} b_k e_k)$
then
$
\sin (\frac{\pi}{2k}) b_k= \frac{1}{k}
$
hence
\[
\lim_{k \so \infty} b_k = \lim_{k \so \infty} \frac{\pi}{2k}
\frac{1}{\sin(\frac{\pi}{2k})} \frac{2}{\pi}=\frac{2}{\pi}.
\]
This means that  $\gamma(\frac{\pi}{2})$ is  epifocal.} 
\end{ex}
This example shows that there exist  epifocal  points which are not 
monofocal. We call them
pathological points. If the exponential map of a Hilbert manifold
has only a finite number of phatological points we will say that the 
exponential map is {\em almost non singular.}   
Clearly, if  the exponential map were Fredholm, 
and this one  must be  of zero index, 
monoconjugate 
points and epiconjugate points along geodesics would coincide. This holds for 
the Hilbert manifold $\Omega (M),$ the free loop space of a compact manifold 
(see \cite{Mi}). Moreover, any finite geodesic in $\Omega (M)$ contains 
at most finitely many points which are conjugate.

Now we shall prove  the Index lemma. This lemma allows us to extend
Rauch and Berger theorems in infinite dimensional context.

Let $X:[0,1] \lra \tpm$, such that $X(0) \in T_p N.$ We define
the {\em focal index} of  $X$ as follows:
$$
\begin{array}{lcl}
I^N (X,X)& =& \int_{0}^{1}
(\langle  \fis X (t), \fis X (t) \rangle - \langle R_t(X(t)),X(t)
\rangle )dt
\\
& - &  \langle \We (X(0)),X(0) \rangle .
\end{array}
$$
If $N=p,$   the Focal Index is called Index and we will 
denote it as $I(X,X).$  We note that any vector fields along
$\gamma$ is a parallel transport of a unique application 
$X:[0,b] \lra \tpm$; we will denote by
$\overline X (t)=\pt{0}{t}(X)$ the vector field along $\gamma$ starting from
$X$.
\begin{lemma}
$I^N (X,X)= D^2 E(\gamma)(\overline X, \overline X)$,  
where  $D^2 E(\gamma)$ is the
index form  of  $B=N \times M \hookrightarrow M \times M.$
\end{lemma}
\n
\dem $\ $ we recall that
\begin{small}
\begin{eqnarray}
\D( \overline X , \overline X ) &=& \int_a^b
\pr \nbp{t}{\overline X (t)} \pr^2
                    \ - \ \langle \overline X (t), 
R(\overline X (t),\fis c(t)) \fis c(t) \rangle dt 
\nonumber \\
             &-& <<A_{(\fis c(a),-\fis c(b))}(\overline X(a), 
\overline X (b)),(\overline X (a),\overline X (b)>>,
             \nonumber
\end{eqnarray}
\end{small}
see \cite{Sa}, where $A$ is the Weingarter operator of 
$N \times M \hookrightarrow M \times M.$
By the above expression, it is enough to prove that
$\nbp{t}{\ov{X}}(t)=\pt{0}{t}(\fis X (t))$. Let $Z(t)$ be a parallel transport
of a vector $Z \in \tpm.$ Then
$$
\begin{array}{ccl}
\langle \nbp{t}{\overline X (t)},Z(t) \rangle  &= & 
\langle \overline X (t),Z(t) \rangle ' \\
                          &= & \langle \fis{X} (t), Z \rangle \\
                          &= & \langle \pt{0}{t}(\fis{X}(t)), Z(t) \rangle.
\ \mathrm{QED}
\end{array}$$

We here give the Focal Index lemma formulated as in the finite 
dimensional Riemannian geometry.
\begin{lemma} \label{fp}
Let $X:[0,b] \lra \tpm$ be a piecewise differential application with
$X(0) \in T_p N$.
Suppose that $T(t)$ is invertible in $(0,a)$. Hence
\[
I^N (X,X) \geq I^N (J,J),
\]
where $J(t)=T(t)(u)$ with $X(b)=T(b)u$. The equality holds if and only if
$X=T(t)(u)$. Hence, if there aren't any focal points along $\gamma$, 
the index of a vector fields $Y$ along $\gamma$ is bigger than the focal index 
of the $N$-Jacobi field $J$ along $\gamma$ such that $W(a)=J(a).$
\end{lemma}
\n
\dem$\ $   $T(t)$ is invertible, then 
there exists a piecewise differential application 
$Y:[0,b] \lra \tpm$ such that
$Y(0)=X(0) \in T_p N$ and 
$X(t)=T(t)(Y(t))$. Hence
\[
\fis{X}(t)=T'(t)((Y(t)) \ + \ T(t)(\fis{Y}(t))=A(t) \ + \ B(t).
\]
The focal index of $X$ is given by
\begin{small}
$$
\begin{array}{ccl}
I(X,X) &=& \int_{0}^{b} ( \langle A(t),A(t) \rangle  \ + \ 2 \langle A(t),B(t) \rangle \ + \
               \langle B(t),B(t) \rangle  dt  \\
       &-&  \int_{0}^{b} \langle R_t (T(t)(Y(t)),T(t)(Y(t)) \rangle dt \ - \
\langle \We (X(0),X(0) \rangle.
\end{array}$$
\end{small}
A straightforward computation show that
\begin{eqnarray}
\langle A(t),A(t) \rangle  &=& \langle T(t)(Y(t)),T'(t)(Y(t) \rangle'   \nonumber  \\
            &-& \langle B(t),A(t) \rangle                 \nonumber   \\
            &+& \langle T(t)(Y(t)),R_t (T(t)(Y(t)) \rangle  \nonumber     \\
            &-& \langle T(t)(\fis{Y}(t)),T'(t)(Y(t)) \rangle
            = ( \langle A(t),B(t) \rangle ),  \nonumber
            \nonumber
\end{eqnarray}
where the last equality depends on the fact that 
$\Phi(t)$ is a family of symmetric bilinear form. Hence, 
the focal index of $X$ is given by
\[
I^N (X,X) \  = \  \langle T(1)(u),T'(1)(u) \rangle \ + \
  \int_{0}^{b} \pr T(t)(\fis{Y}(t)) \pr^{2}dt .
\]
thus proving our lemma. QED 
\begin{cor} \label{compara}
Let $(M,g)$ be a Riemannian manifold and let $S$ and $\Sigma$ be submanifolds
of codimension $1.$ We denote by $N$ and $\ov{N}$ the normal vector fields
respectively in $S$ and $\Sigma$. Suppose that  in some point 
$p \in S \cap \Sigma$ we have    $N_p=\ov{N}_p$ and 
\[
g(\nb_{X}{N},X)  < g(\nb_{X}{\ov{N}},X),
\]
for every  $X \in T_p \Sigma =T_p S$. Then, if the Jacobi flow $T$ of
$S$ is invertible in  $(0,b)$ then the Jacobi flow of
$\Sigma$ must be injective in $(0,b).$ Moreover, if 
$A-\ov{A}$ is invertible, where 
$A$ and $\ov{A}$ are the Weingarter operators in
$p$ respectively of $S$ and $\Sigma,$ then the Jacobi flow of $\Sigma$ is also
invertible in $(0,b).$ 
\end{cor}
\n
\dem $\ $ let $Y(t)$ be $\Sigma$-Jacobi field. 
Since $T(t)$ is invertible then there exists a piecewise
application 
$X:[0,s] \lra T_p M$  with 
$X(0) \in T_p S$ such that
$Y(t)=T(t)(X(t)).$ Hence
$$
\begin{array}{rcl}
Y(0) &=& T(0)(X(0)) \\
\fis Y (0)&=& T'(0)(X(0))\ + \ T(0)(\fis X (0)) \\
(-\ov{A}(Y(0)), \fis Y (0)^n) &=& (-A (X(0))\ +\ \fis X (0)^t,0)).
\end{array}
$$
Hence, $Y(0)=X(0)$ and the tangent component 
of $\fis X (0)$ is given by $(A\ -\ \ov{A})(X(0)).$
Then, 
\begin{small}
$$
\begin{array}{lcl}
g(Y(s),\nbp{s}{Y} (s))& = & I^S (Y,Y) \\
                     & = & g((A -  \ov{A})(X(0)), X(0)) +
                        \int_{0}^{s} \pr T(t)(\fis X (t)) \pr^2 dt  \\
                     & > & 0.
\end{array}
$$
\end{small}
In particular, the Jacobi flow of $\Sigma$ is injective in $(0,b).$
Moreover, if $A- \ov{A}$ is invertible, using the above formula, 
it is easy to prove that the 
image of the Jacobi flow relative to $\Sigma$ is a closed subspace for
every $t\in (0,b)$ and using $(i)$ of proposition \ref{focal}
we get our claim. QED \\

Now, assume that there exists a pathological point on the interior of 
$\gamma;$ this means that the Jacobi flow is an isomorphism
for every $ t \neq t_o$ in  $(0,b)$ and
in $t_o,$ $T(t_o)$ is a linear operator whose image is a dense subspace.
Let $X:[0,b] \lra \tpm$ be a piecewise application with
$X(0) \in T_p N$. 
Given $\epsilon>0$, there exist  $X_n^{\epsilon}$, $n=1,2$
such that
\[
\begin{array}{ccc}
& \pr T(\to)( X_1^{\epsilon})) -  X(\to)) \pr   \leq \frac{\epsilon}{4} & \\
& \pr T(\to)( X_2^{\epsilon})) -  \fis X(\to)) \pr \leq \frac{\epsilon}{4}.
\end{array}
\]
Choose  $Y^{\epsilon}$ such that
\[
\pr T(\to)(Y^{\epsilon}) -  T'(\to)(X_1^{\epsilon}) \pr \leq
\frac{\epsilon}{4} .
\]
Hence there exists  $\eta(\epsilon) \leq \frac{\epsilon}{2}$
such that for $ t \in (\eta(\epsilon) -  \to, \eta(\epsilon) +  \to)$ 
we have
\[
\begin{array}{ccc}
(1) & \pr T(t)( X_1^{\epsilon}  + 
   (t  -  \to)(X_2^{\epsilon}   -  Y^{\epsilon}))
    -   X(t)  \pr \leq \epsilon, & \\
(2) & \pr \d1{}{t} (T(t)( X_1^{\epsilon}  + 
(t  -  \to)(X_2^{\epsilon}  - Y^{\epsilon}))) 
  -  \fis X (t) \pr \leq \epsilon. &
\end{array}
\]
We denote by  $X^{\epsilon}$ the application
\begin{small}
\[
X^{\epsilon}(t)=
\left \{ \begin{array}{ll}
X(t) & {\rm if} \ \ 0 \leq t \leq \to\ - \
\eta( \epsilon); \\
T(t)(X_1^{\epsilon}  + 
(t  -  \to)(X_2^{\epsilon}   - Y^{\epsilon})) & {\rm if} \
\to  - \eta( \epsilon) < t < \to  +  \eta( \epsilon); \\
X(t)   & {\rm if}\  \to -  \eta( \epsilon ) \leq t \leq b .
\end{array}
\right. \]
\end{small}
Clearly,  $X^{\epsilon}=T(t)(Y(t)),$ where  $Y(t)$ is again a 
piecewise application,
at less at the points $t=\to +  (\eta (\epsilon)$ and 
$t=\to - \eta(\epsilon),$ 
and using the same arguments in  lemma \ref{fp}, we shall prove
\[
I^{N}(X^{\epsilon},X^{\epsilon}) \geq I^{N}(T(t)(u),T(t)(u)) .
\]
On the other hand, the Focal Index of $X$ is given by
\begin{small}
$$\begin{array}{lcl}
I(X,X)&=& I(X^{\epsilon},X^{\epsilon}) \\
      &-& \int_{ \to - \eta(\epsilon) }^{ \to + \eta(\epsilon) }
          \langle \fis{X}^{\epsilon}(t),\fis{X}^{\epsilon}(t) \rangle\ - \
 \langle R({X}^{\epsilon}(t), \fis c (t))\fis c (t),{X}^{\epsilon}(t) \rangle dt \\
      &+&  \int_{ \to - \eta ( \epsilon)}^{ \to + \eta(\epsilon) }
          \langle \fis X  (t), \fis X (t) \rangle \ - \
          \langle R(X(t), \fis c (t)) \fis c (t),X(t) \rangle dt .
\end{array}$$
\end{small}
Now, using $(1)$ and $(2)$  it is easy to check that
$$
\lim_{\epsilon \so 0} I^N ( X^{\epsilon}, X^{\epsilon} )= I^N (X,X)
\geq I^N (J,J),
$$ 
where $J(t)=T(t)(u).$ This proves the Focal Index lemma when there is  
a pathological point. Clearly, the same proof can be generalized to a finite 
number of pathological points.
\begin{lemma} {\bf (Focal Index lemma)}
Let $\gamma:[0,b] \lra M$ be a geodesic with a finite number of pathological
points. Then for every vector fields $X$ along  $\gamma$ with
$X(0) \in T_p N,$ the index form of $X$ relative to  the submanifold
$N \times M \hookrightarrow M \times M$,
satisfies $\D (X,X) \geq \D(J,J)$, where
$J$  is the $N$-Jacobi field such that  $J(b)=X(b)$.
\end{lemma}
In particular, if $\gamma$ has a finite number of pathological
points, then $\gamma$ is a local minimum and the Index form
$\D$ is non negative defined. 
\begin{cor}
Let $\gamma:[0, \infty) \lra M$ be a geodesic, and let $\gamma (t_o)$ 
be a monoconjugate point. Then $\gamma:[0,t] \lra M$ is not minimal for
$t>t_o.$
\end{cor}
Now, we shall prove the Rauch and Berger 
theorem and several corollaries. However, the proof are almost the same
as the finite dimensional case, since these can be proved 
using Focal Index lemma: 
then we will give only a brief proof of the Rauch theorem. 
\begin{thm}{\bf (Rauch)}\label{R}
Let $(M,\langle \ ,\ \rangle)$, $(N,\langle \ ,\ \rangle^*)$
be Hilbert manifolds
modeled on $\H_1$ and $\H_2$ respectively, with $\H_1$
isometric to a closed  subspace of $\H_2$. Let
\[
c:[0,a] \lra M,\ \ c^*:[0,a] \lra N
\]
be two geodesics of same length.
We assume that $c^*$ has at most a finite number of pathological 
points in its interior. Suppose also that for every  $t\in [0,a]$
and for every  $X \in T_{c(t)}M$, $X_o \in T_{c^*(t)}N$ we have
\[
K^N (X_o,\fis c^*(t)) \geq K^M (X,\fis c(t)).
\]
Let  $J$ and  $J^*$ be Jacobi fields along $c$ and $c^*$ such that
$J(0)$ and $J^*(0)$ are tangent to $c$ and $c^*$ respectively and 
\begin{enumerate}
\item $\pr J(0) \pr = \pr J^*(0) \pr$;
\item $\langle \fis c(0),\nbp {t}{J(0)} \rangle =
\langle \fis c^*(0),\nbp {t}{J^*(0)} \rangle $;
\item $\pr \nbp {t}{J(0)} \pr =\pr \nbp {t}{J^*(0)} \pr$.
\end{enumerate}
Hence, for every  $t\in [0,a]$
\[
\pr J(t) \pr \geq \pr J^*(t) \pr.
\]
\end{thm}
\n
\dem $\ $ 
It is easy to check that we will restrict ourself to the case in which 
the  Jacobi fields satisfy the following condition:
\[
0=\pr J(0) \pr=\langle \fis c(0),\nbp {t}{J(0)} \rangle =  \pr J^*(0) 
\pr= \langle \fis
c(0),\nbp{t}{J^*(0)}.  \rangle
\]
Note, by assumption, $J^*(t) \neq 0$. Let
$\to \in [0,a]$ be an isometry
$$
\begin{array}{lcl}
&F: T_{c(0)}M \lra T_{c^*(0)}N &  \\
&F(\fis c(0))=\fis {c}^*(0) &  \\
&F(\pt{\to}{0}(J(\to))=\ptn{\to}{0}(J^*(\to))
\frac{\pr J(\to) \pr}{\pr J^*(\to) \pr}.
\end{array}
$$
We denote by
\begin{eqnarray}
& i_t: T_{c(t)}M \lra T_{c^*(t)}N \nonumber \\
& i_t=\ptn {0}{t} \circ F \circ \pt {t}{0}, \nonumber
\end{eqnarray}
a family of isometries for $0 \leq t \leq t_o;$  
it 
is easy to check that $i_t$ commute with the Levi Civita connection. 
Let $W(t)=i_t(J(t)).$ Then
\begin{small}
$$\begin{array}{lcl}
D^2 E(c^*)(W,W) &=& \int_0^{\to} \pr \nbp{t}{W(t)} \pr^2 \
                  - \ \langle R^N (\fis c^*(t),W(t))\fis c^*(t),W \rangle dt  \\
              &\leq & \int_0^{\to} \pr \nbp {t}{J(t)} \pr^2
              - \langle R^M (\fis c(t),J(t))\fis c(t),J(t) \rangle dt  \\
            &=& \D(J,J) .
\end{array}
$$
\end{small}
On the other hand,
\begin{eqnarray}
\frac{1}{2}\d1{}{t} \mid_{t=\to} \langle J(t),J(t) \rangle &=& \langle J(\to),\nbp{t}{J(\to)}
\rangle
 \nonumber  \\
         &=& \D(J,J) \nonumber \\
  &\geq& D^2 E(c^*)(W,W) \nonumber \\
  & \geq  & 
   \d1{}{t} \mid_{t=\to} \langle J^*(t),J^*(t) \rangle^*
         \frac{\pr J(\to) \pr^2 }{\pr J^*(\to) \pr^2 } . \nonumber
\end{eqnarray}
where the last before inequality follows from the Focal Index lemma.
Now, let $\epsilon >0$. For every $t \geq \epsilon$ we have
\[
\d1{}{t} \log (\pr J(t)\pr^2) \geq \d1{}{t} \log (\pr J^*(t) \pr^2),
\]
that implies
\[
\frac{\pr J(t) \pr^2 }{\pr J(\epsilon) \pr^2 } \geq
\frac{\pr J^*(t) \pr^2 }{\pr J^*(\epsilon) \pr^2} .
\]
By $\pr \nbp {t}{J(0)} \pr = \pr \nbp {t}{J^*(0)} \pr,$ 
when $\epsilon \so 0$ we get our claim.
QED\\

\n
Misiolek, see \cite{Mi}, proved that in $\Omega(M)$ the index of any
finite geodesic if finite. By Rauch theorem, we have that the 
sectional curvature of $\Omega (M),$ with the $H^1$ metric, 
cannot be positive along any geodesic. 
Indeed, if $K \geq K_o>0,$ we are able to compare $\Omega(M)$  with 
the sphere of radius $\frac{1}{\sqrt{K_o}}$ and we may prove that 
the index along any geodesic of length bigger than
$t=\frac{\pi}{\sqrt{K_o}}$ is infinite. 
\begin{cor}\label{cr}
Let  $M$, $N$ be Hilbert manifolds modeled on $\H_1$ and  $\H_2$
where $\H_2$ is isometric to a closed subspace of $\H_1$. Assume that for
every
$m \in M$ and $n \in N$ and for every $\eta \in  \tpm$  e
$\beta \in T_{n}N$ 2-subspaces we have
\[
K^M (\eta) \geq K^N (\beta) .
\]
Let  $i:T_n N \lra \tpm$ be an isometry and let $r>0$ such that
\begin{center}
\begin{tabular}{lcl}
&$\e n:\bn{n}{r} \lra \bd{B}{m}{r}$ & is  a diffeomorphism \\
&$\e m:\bn{m}{r} \lra \bd{B}{n}{r}$ & is almost non singular. \\
\end{tabular}
\end{center}
Let
Let $c:[a,b] \lra  \bn{n}{r}$ be a piecewise curve. Then
\[
L(\e n (c)) \geq L(\e m (i\circ c)) .
\]
\end{cor}
\begin{thm}\label{B}{\bf (Berger)}
Let $(M,g)$ and $(N,h)$ be an Hilbert manifolds modeled on 
$\H_1$ e $\H_2$, where $\H_1$ is isometric to a closed subspace of
$\H_2$. Let $\gamma_1:[0,b] \lra M$ and $\gamma_2:[0,b] \lra N$ be two 
geodesics with the same length. Assume that for every
$X_1 \in T_{\gamma_1(t)}M$ and
$X_2 \in T_{\gamma_2(t)}N$ we have
\[
K^N(X_2,\fis{\gamma_1}(t)) \geq K^M (X_1,\fis{\gamma_2}(t)),\ 
\langle X_1,\fis{\gamma_1}(t))
\rangle = \langle X_2,\fis{\gamma_2}(t) \rangle =0.
\]
Assume furthermore that $\gamma_2$ has at most a finite number of
pathological points, on its interior,
of the  geodesic submanifold $N$ defined by $\fis{\gamma_2}(0)$.
Let $J$ and
$J^*$ Jacobi fields along $\gamma_1$ and $\gamma_2$ 
satisfying  $\nbp{t}{J(0)}$  and $\nbp{t}{J^* (0)}$ are tangent to $\ga 1$ and
$\ga 2$ and 
\begin{enumerate}
\item $\pr \nbp{t}{J(0)} \pr = \pr \nbp{t}{J^* (0)} \pr$,
\item $\langle \fis{\gamma_1}(0),J(0) \rangle = \langle 
\fis{\gamma_2}(0),J^*(0) \rangle ,\ \pr J(0) \pr = \pr J^*(0) \pr$.
\end{enumerate}
Then
\[
\pr J(t) \pr \geq \pr J^* (t) \pr, 
\]
for every $t \in [0,b]$.
\end{thm}
\begin{cor}\label{cb}
Let $(M,g)$ and $(N,h)$ be an Hilbert manifolds modeled on 
$\H_1$ and $\H_2$ respectively , where $\H_1$ is isometric to a 
closed subspace of
$\H_2$. Let $\gamma_1:[0,b] \lra M$ and $\gamma_2:[0,b] \lra N$ be two 
geodesics with the same length. Let $V_1(t)$ and  $V_2(t)$ be parallel
unit vectors  along  $\gamma_1$ and $\gamma_2$ 
which are everywhere perpendicular to the tangent vectors of $\gamma_1$ and
$\gamma_2.$ Let $f:I \lra \R$ be a positive function and let
\begin{center}
$b(t)=\e{\gamma_1(t)} (f(t)V_1(t))$,\\
$b^*(t)=\e{\gamma_2(t)}(f(t)V_2(t))$,
\end{center}
two curves. Assume that $K^N \geq K^M$ and for any $t \in I$ the geodesics 
\[
\eta_o (s)=\e{\gamma_2(t)}(sf(t)V_2(t)),\ 0 \leq s \leq 1
\]
contains no focal points of the geodesic submanifold defined by
$\fis{\eta} (0).$ Then $L(b) \geq L(b^*)$.
\end{cor}
\begin{cor}\label{zio}
Let $(M,g)$ be an Hilbert manifold such that $H \leq K \leq L,$ $H>0$ 
and let $\gamma:[0,b] \lra M$ be a unit speed geodesic.
Then
\begin{enumerate}
\item the distance $d,$ along $\gamma,$ from $\gamma (0)$ 
to the first monoconjugate or epiconjugate 
satisfies the following inequality
$$
\frac{\pi}{\sqrt{H}} \leq d \leq \frac{\pi}{\sqrt{L}};
$$
\item the distance $d$, along $\gamma,$ from $\gamma (0)$ to 
first monofocal or epifocal point, of the
geodesic submanifold defined by $\fis{\gamma} (0),$ 
satisfies the following inequality
$$
\frac{\pi}{2 \sqrt{H}} \leq d \leq \frac{\pi}{2 \sqrt{L}}.
$$
\end{enumerate}
\end{cor}   
\section{Hilbert Manifolds: a global theory}
\mbox{}
The Rauch and Berger theorems are very important to understand 
the geometry of the complete manifolds with upper or 
lower curvature bounded. Indeed, 
we can compare these manifolds with the complete Hilbert manifolds with 
constant 
curvature and the geometry  of these is well known. We saw that in a 
complete Hilbert manifold the exponential map may not be surjective.  
When the curvature is upper bounded by a constant we have the following
result.  
\begin{prop} \label{s1}
Let $(M,g)$ be a complete Hilbert manifold such that $K \leq H$. If
$c:[0,1] \lra M$ is a piecewise differential curve, with
$L(c)< \frac{\pi}{\sqrt{H}}$ if $H>0$,
then there exists a unique piecewise differential  curve 
$\overline{c}:[0,1] \lra \overline{\bn{c(0)}{L(c)}}$
such that $\e {c(0)} ( \overline{c(t)})=c(t)$. In particular,
\[
\e p(\bn{p}{r})=\bd{\cal B}{p}{r}
\]
for every $p \in M$ and $r \geq 0$ if $H \leq 0$ or
$r< \frac{\pi}{\sqrt{H}}$ if $H>0$.
\end{prop}
\n
\dem $\ $ we will give the proof only when $H>0;$ the other case is similar. 
Take 
\begin{small}
\[
\to=sup \{ \ t \in [0,1]:\ \exists! \  \overline{c}:[0,t] \lra
\overline{\bn{c(0)}{L(c)}}, \  {\mathrm with}\  \e p 
( \overline{c} (t) )=c(t) \};
\]
\end{small}
\n
$t_o$ is positive by Rauch theorem, and we shall prove that $t_o$ is in fact 
$1.$ Let
\[
\overline{c}:[0,\to) \lra \overline{ \bn{c(0)}{L(c)} }
\]
the unique lift of $c$; using Rauch weak theorem we have
$$
\begin{array}{rcl}
\pr \fis c (t) \pr &=& \pr d ( \e{c(o)})_{\overline{c(t)}}
(\fis{\overline {c}}(t)) \pr \\
&\geq & \frac {\sin(\pr \overline{c(t)})
 \pr \sqrt{H} ) }{\sqrt{H} \pr \overline{c(t)} \pr }
 \pr \fis{\overline{c}}(t) \pr \\
&\geq & \frac{\sin(L(c))\sqrt{H})}{\sqrt{H}L(c)}
\pr \fis{\overline{c}}(t) \pr,
\end{array}$$
so we get
\[
\lim_{t \so \to} \int_{0}^{\to} \pr \fis{\overline{c}}(t) \pr dt \
<  \infty.
\]
However, $\overline{\bn{c(0)}{L(c)}}$ is a complete metric space so
$\lim_{t \so \to} \overline{c}(t)=q$ and by Gauss lemma and the definition of
$t_o$ we get  $t_o=1.$ 
QED
\begin{cor}
Let $(M,g)$ be a complete Hilbert manifold such $K \leq H.$ Let
$p,q \in M$, such that  $d(p,q)< \frac{\pi}{\sqrt{H}}$ if $H>0$.
Hence, at least one of the following facts holds:
\begin{enumerate}
\item  there exists a unique minimal geodesic between $p,q$;
\item  there exists a sequence  $\gamma_n$ of geodesics
 from $p$ to $q$  such that $L(\gamma_n)>L(\gamma_{n+1})$ and
$L(\gamma_n) \so d(p,q)$.
\end{enumerate}
\end{cor}
Next we claim a very useful lemma that we will use in the following
proofs.
\begin{lemma} \label{s2}
Let $(M,g)$ be an Hilbert manifold such that $K \geq L>0.$
Suppose there exists  a point $p \in M$ on which $\e p$
 is almost non singular in
$\bn{p}{r}$. Let $\delta(s)$ be a curve joining two antipodal points
on the sphere of radius $s$ in $\tpm$. Let
$\Delta$ be
the Image, via $\e p,$ of the curve $\delta(s).$ Then
\[
L(\Delta)\leq \frac{\pi}{\sqrt{L}} \sin (s \sqrt{L} ),
\]
for $s< r$.
\end{lemma}
\n
\dem $\ $ Let $S_{\frac{1}{\sqrt{L}}} (\tpm \times \R)$ be the sphere of radius
$\frac{1}{\sqrt{L}}$  and let
$N=(0,\frac{1}{\sqrt{L}}) \in S_{\frac{1}{\sqrt{L}}} (\tpm \times \R)$;
it is easy to check that
\[
\e N (v)=
\frac{1}{\sqrt{L}}(\cos ( \pr v \pr \sqrt{L} )N\ +\
\sin ( \pr v \pr  \sqrt{L})
\frac{v}{\pr v \pr}).
\]
Let $v,w \in \tpm$ be such that $\langle v,w \rangle =0,
\ \langle v,v \rangle = \langle w,w \rangle
=1$;
any meridian on the sphere of radius $s$ can be parametrized as follows:
\[
c_s(t)=s( \cos (s)v \ + \ \sin (s)w )
\]
and $L(\e N (c_s))= \frac{\pi}{\sqrt{L}} \sin (s \sqrt{L})$. By
corollary \ref{cr}, we have 
\[
L(\Delta) \leq \frac{\pi}{\sqrt{L}} \sin (s \sqrt{L}).
\]
QED 
\begin{prop} \label{bt1}
Let $(M,g)$ be an Hilbert manifold such that $K \geq 1$. Suppose there exists
a point $p \in M$ such that $\e p$ is almost non singular in 
$\bn{p}{\pi}.$ Then $M$ has constant curvature $K=1$ and is covered by the 
unit sphere $S(\tpm \times \R)$. Furthermore, $M$ is a complete Hilbert
manifold.
\end{prop}
\n
\dem $\ $
let $S_r(\tpm),$  $ r < \pi$ be the sphere of radius $r$ in $\tpm$; 
using lemma \ref{s2} we get that 
the diameter of $\e p ( S_r(\tpm) ) \so 0$ when $r \so \pi.$ Hence
$\e p ( S_{\pi} (\tpm))=q.$
Let $N=(0,1) \in  S(\tpm \times \R).$ We define
\[
\phi(m)= \left \{ \begin{array}{l}
\e p ({\e N}^{-1}(m)) \ \           {\rm se} \  m \neq -N; \\
q \ \ \ \ \ \ \ \ \ \ \ \ \ \ \ \ \ \ \ \ \  \ \ \   {\rm se} \ m=-N.
\end{array}
\right.
\]
We claim that  $\phi$ is a local isometry. Firstly, note that 
any geodesic which starts at $p$ get to $q.$ Hence any Jacobi field along a 
geodesic which starts in $p$ is zero in $q.$ Moreover the Index form
of any geodesic $\gamma: [0,t] \lra M,$ $t \leq \pi$ and $\gamma (0)=p,$ 
is non negative definite because $\e p$ is almost non singular in 
$\bn{p}{\pi}.$
Let $\gamma :[0,\pi] \lra M$
be a geodesic such that $ \gamma (0)=p.$ Let $W(t)$
be a parallel transport along $\gamma$ of a unitary and perpendicular vector
to $\fis{\gamma}(0)$. The index form of $W$ along $\gamma$ is given by 
\begin{small}
$$
\begin{array}{ccl}
0  & \leq & \D (Y,Y) \\
& =    & \int_{0}^{\pi} \langle Y(t),Y(t) \rangle  \ - \
            \langle R(Y(t),\fis{\gamma}(t))\fis{\gamma}(t),Y(t) \rangle dt \\
   & \leq & \int_{0}^{\pi} ( \cos^2 t\ - \ \sin^2 t ) dt  \\
   & =    & 0 .
\end{array}$$
\end{small}
Hence $K(W(t),\fis{c}(t))=1$ and  $Y(t)$ is a Jabobi field. Now, using Cartan
theorem, about local isometry, and  proposition 6.9 pag. 222 in \cite{La}, 
we get our result.
QED\\

Next, we claim the Berger-Topogonov theorem of maximal diameter. 
\begin{thm} \label{bt} {\bf (Berger-Topogonov)}
Let $(M,g)$ be a complete Hilbert manifold such that 
$\delta \leq K \leq 1$. Suppose that there exist
two points $p,q$ with $d(p,q)=\frac{\pi}{\sqrt{\delta}}$ and at least a 
minimal geodesic from $p$ to $q.$ Then $M$ is
isometric to a sphere $S_{\frac{1}{\sqrt{\delta}}}(\tpm \times \R)$.
\end{thm}
\n
\dem $\ $
By corollary  \ref{zio}, the distance from the first focal point 
along any geodesic 
is at
least $\frac{\pi}{2};$ furthermore 
by Bonnet theorem, $d(M)=\frac{\pi}{\sqrt{\delta}}$.  Let
$\gamma:[0,\frac{\pi}{\sqrt{\delta}}] \lra M$ be a minimal geodesic from
$p$ to $q$. Take the following vector field along $\gamma$
\[
Y(t)=\frac{1}{\sqrt{\delta}} \sin (t \sqrt{\delta}) W(t),
\]
where $W(t)$ is the parallel transport along $\gamma$ perpendicular to
$\fis{\gamma}(t)$. We define the following variation of $\gamma$ to be
\[
\Omega (s,t)=\e{\gamma(t)}(sY(t)),\ 0 \leq s \leq 1,\  0 \leq t \leq
\frac{\pi}{\sqrt{\delta}}.
\]
Any curve $\Omega (s,\cdot)$ joins $p,q$ and by corollary \ref{cr} we get
\[
L[\Omega (s,\cdot )] \leq  \frac{\pi}{\sqrt{\delta}}
\]
Hence  any curve is a minimal geodesic and $Y(t)$ is a Jacobi field. 
Furthermore $\Omega (\cdot,\cdot)$
is a totally geodesic submanifold. 
Now, it easy to check that $\e p$ is non singular and injective in 
$\bn{p}{\frac{\pi}{\sqrt{\delta}}}$, 
$\e p (S_{\frac{\pi}{\sqrt{\delta}}} (\tpm) )=q$ and the application
\[
\Phi:S_{\frac{1}{\sqrt{\delta}}} (\tpm \times \R) \lra M
\]
defined by
\[
\left \{ \begin{array}{ll}
\Phi(m)= \e p \circ {\e N}^{-1}(m) & {\rm se} \  m \neq -N ; \\
\Phi(-N)= q                    & {\rm se}  \ m=-N;
\end{array}
\right.
\]
where $N=(0,\frac{1}{\sqrt{\delta}})$, is an isometry. QED \\

The Sphere theorem is one of the most beautiful theorem in classic Riemannian
geometry. Unfortunately, we haven't found yet a proof in infinite dimensional
case. We will show two  soft versions of Sphere theorem: one is the Sphere 
theorem due in finite dimensional Riemannian geometry by Rauch 
(see \cite{Ra}), 
with the strong assumption on the injectivity 
radius, with pinching $\sim \frac{3}{4}$ and the other is the Sphere
theorem in the class of Hopf-Rinow manifolds, on which we shall prove that 
the
Topogonov theorem holds, with pinching $\frac{4}{9}.$ In the first case 
the fundamental lemma is the following.  
\begin{lemma}
Let $(M,g)$ be an Hilbert manifold such that
$K \leq H$, $H>0.$
Let  $\phi: S(\H) \lra M$, where $S(\H)$ is an unit sphere in an 
Hilbert space, 
be a local homeomorphism onto the image such that:
$\phi(N)=p$ and the image of every meridian is a curve of length
$r\leq r_o < \frac{\pi}{\sqrt{H}}$ .
Then there exists a locally homeomorphism 
\[
\ov{\phi}: S(\H) \lra \ov{ \bn{p}{r_o} },
\]
such that $\e p \circ \ov{\phi} = \phi$.
\end{lemma}
\n
\dem $\ $
We apply the proposition \ref{s1} to each meridian and we get
\[
\ov{\phi}: S(\H) - \{ -N \} \lra \ov{ \bn{p}{r_o} },
\]
with $\e p \circ \ov{\phi} = \phi.$ We claim that 
we can extend $\ov{\phi}$ to $-N.$ 
Let $\xi(t)=\e N (tv) $ be a meridian starting from 
$N.$ Let 
$\gamma(t)=\phi(\xi(t))$ and let $\ov{\gamma}(t)$ the lift of $\gamma$. By
assumption, for every $t \in [0,\pi]$ there exists an open subset
$W(t)$ of $\ov{\gamma}(t)$ and an open subset $U(t)$ of
$\gamma(t)$ such that
\[
\e p : W(t) \lra U(t)
\]
is an onto diffeomorphism. Now, $\phi$  is a local homeorphism then
there exists an open subset $V(t)$ of $\xi (t)$ 
such that  $\phi(V(t)) \subset U(t)$ and 
$\phi$ on $V(t)$ is an homeomorphism. The closed interval is compact, then
there exits a partition
$0=t_o \leq t_1 \leq  \ldots \leq t_{n-1} \leq t_{n}=\pi$ such that
\[
\left \{
\begin{array}{l}
\xi( [0,\pi]) \subseteq V(t_o) \cup \cdots \cup V(t_n)=V; \\
\gamma( [0,\pi]) \subseteq U(t_o) \cup \cdots \cup U(t_n)=U; \\
\ov{\gamma}([0, \pi]) \subseteq W(t_o) \cup \cdots \cup W(t_n)=W; \\
\end{array}
\right.
\]
and another partition
$ 0 < s_1 < \cdots < s_{n} < s_{n+1}=\pi$
such that
$ t_i < s_{i+1} < t_{i+1}$ and $\xi(s_{i+1}) \in V(t_i) \cap V( t_{i+1})$ for
$ 0 \leq i \leq n-1$. Now it is easy to see that
\[
C:= \{ w \in T_N S(\H) :\ \ov{\e N (tw)}(\pi)=\ov{\xi}(\pi) \}
\]
where $\ov{\xi}$ is the unique lift of $\xi,$ is open and closed. Hence
$\ov{\phi}$ can be extended in $-N$ and $\phi$ is a local  homeomorphism.
QED \\

Now we claim that a manifold with pinching $\sim \frac{3}{4}$ 
can be covered by two geodesic balls.

Let $(M,g)$ be a complete Hilbert manifold. Suppose that the sectional
curvature satisfies the inequality 
$0<h \leq K \leq 1,$  where $h$ is a solution of
the equation
\[
\sin \sqrt{\pi h}=\frac{\sqrt{h}}{2} \ \ (h \sim \frac{3}{4}).
\]
Let $p \in M$. By  Rauch  theorem we get that on the geodesic ball
of radius  $\pi$ there aren't conjugate points. We denote by 
$\Delta$,
the meridian of the sphere in $\tpm$ of radius $\pi.$ Then 
\[
L[\Delta] \leq \frac{\pi}{\sqrt{h}}\sin {\pi \sqrt{h}} \leq
\frac{\pi}{2} .
\]
In particular, there exists $\epsilon>0$ such that the image of any
meridian in the sphere of radius $\pi - \epsilon$ is a curve with length
$r\leq r_1 < r_o < \pi-\epsilon$. Furthermore, $\epsilon$ 
does not  depend from $p.$
\begin{lemma}
Let $p \in M$  and let  $q \in \e p (S_{\pi- \epsilon} (\tpm))$. Then
\[
M=\ov{\bd{\cal B}{p}{\pi-\epsilon}} \bigcup
\ov{\bd{\cal B}{q}{\pi - \epsilon}}.
\]
\end{lemma}
\n
\dem $\ $
take $m \in M$ and let $c:[0,1] \lra M$ be a piecewise curve joining
$p$ and $m$.
Take
\[
t_o=sup \{ t\in [0,1]:\ \exists \ \ov{c}:[0,t] \lra \ov{\bn{p}{\pi- \epsilon}}
:\ \e p (\ov{c}(s))=c(s) \}.
\]
As in proposition \ref{s1}, $\ov{c}$ is defined in
$t_o$ and $L(c_{[0,t_o]}) \geq \pi -  \epsilon$. If $t_o=1$  we get our 
claim; otherwise $c(t_o) \in \e q (\ov{\bn{q}{r_1}}).$ Now, we define
\[
t_1= sup \{ t \geq  t_o: \exists \
\ov{c}:[t_o, t] \lra \ov{\bn{q}{\pi-\epsilon}} .
\}
\]
Now, $r_1< \pi \ - \ \epsilon$ so $t_1>0$  and, as before,
$\ov{c} (t_1)$ is well-defined. If $t_1=1$ we get our result. Otherwise
$\ov{c}(t_1) \in \e p (\ov{\bn{p}{r_1}})$ and by Gauss lemma we get
\[
L[c_{[t_o,t_1]} ] \geq r_1\ - \ r_o.
\]
However, the curve $c$ has finite length, then after a finite numbers of steps
we get that $m$ either belongs to $\ov{\bd{\cal B}{p}{\pi-\epsilon}}$ or
to $\ov{\bd{\cal B}{q}{\pi - \epsilon}}.$ QED \\

Before proving the Sphere Rauch theorem , we recall that the injectivity 
radius of a complete Hilbert manifold is defined by
\begin{small}
$$i(M)= sup \{ r>0: \e p: \bn{p}{r} \lra  
\bd{\cal B}{p}{r} \ \mathrm{is \ a \ diff.\ onto} \ \forall p \in M \}.$$
\end{small}
\begin{thm} {\bf (Sphere Rauch theorem )} Let $(M,g)$ be an Hilbert 
complete manifold
modeled on $\H$ such that $0< h \leq K \leq 1$,
where $h$ is the  solution of the equation 
$\sin (\pi \sqrt{h})=\frac{\sqrt{h}}{2}$.
Assume also that the injectivity radius $i(M) \geq \pi.$ Then 
$M$ is contractible. Furthermore, if $\H$  is a separable Hilbert space 
then $M$ is diffeomorphic  to $S(l_2 )$.
\end{thm}
\n
\dem $\ $
we recall that an infinite dimensional sphere is a deformation retract of
the unit closed disk, because by Bessega theorem, see \cite{BE}, there
exists a diffeomorphism from 
$\H$ to $\H- \{0 \}$ which is the identity outside the unit disk.
When a infinite dimensional manifold $M$ is modeled on Banach space, 
$M$ is contractible if and only if  $\pi_k (M)=0,$ 
for every $k \in \N$ (see \cite{Pa2}).
Let $f:S^k \lra M$ be a continuous application and let
\[
H:\ov{ \bd{\cal B}{p}{ \pi  - \frac{\epsilon}{2} } }  \times [0.1] \lra
  \ov{ \bd{\cal B}{p}{ \pi  - \frac{\epsilon}{2}}}
\]
be the homotopy from the identity map and the retraction on the boundary.
We can extend a map on $M,$ that we denote by $\tilde{H}$,  
fixing the complementary of 
$\bd{\cal B}{p}{\pi - \frac{\epsilon}{2} }$. Then
$$
\begin{array}{lll}
& F: S^k \times [0,1] \lra M,& \\
& F(x,t)=\tilde{H}(f(x),t). & \\
\end{array}
$$
is a homotopy between $f$ and an application
$\tilde{f}:S^k \lra \ov{ \bd{\cal B}{q}{\pi - \epsilon}}$. Then $f$ 
is nullhomotopic.
If $\H$ is separable, using the 
Kuiper-Burghelea theorem (see \cite{KB}), homotopy classifies the Hilbert
manifolds, up to a diffeomorphism, so $M$ is diffeomorphic to the sphere. 
QED \\

Next we claim another version of the Sphere theorem in the class
of Hopf-Rinow manifolds. However, the main result in this class of Hilbert 
manifolds is the Topogonov theorem.  
From our results appeared in section $4,$ we shall prove it using the same 
idea in  \cite{CE} page 42. First of all, 
we start with the following result that we may prove as
in \cite{Kl}  2.7.11 Proposition page 224.   
\begin{lemma} \label{toro}
Let $(M,g)$ be a Hilbert manifold with bounded sectional curvature 
$H \leq K \leq \Delta,$ where $H, \Delta$ are constant. 
Let $\gga$ be a geodesic segment in $M$  such that
$\ga 1 (0) = \ga 2 (0)$ and $\alpha=(- \fis{\ga 1} (l_1), \fis{\ga 2} (0)).$
We call such configuration hinge. Suppose that
$\ga 1$ and  $\ga 2$
are minimal geodesic with perimeter
$P=l_1 \ +\ l_2 \ + \ d(\ga 1 (0),\ga 2 (\lo 2)) \leq
\frac{2 \pi}{\sqrt{H}} \ -  \ 4 \epsilon$,
 $\epsilon>0$ if $H>0.$ In addition,
\begin{itemize}
\item [(i)] if $H \leq 0$ then
$l_2 \leq \frac{\pi}{2 \sqrt{\Delta}}$;
\item [(ii)] if $H>0$ then
\[
l_2 \leq
inf ( \epsilon, \frac{\sin \sqrt {H}\epsilon }{\sqrt{H}}
                  \sin \frac{{\pi} \sqrt{H}}{2 \sqrt{\Delta}},
       \frac{\pi}{2 \sqrt{\Delta}} ).
\]
\end{itemize}
Let $\ggao$ be a hinge $\mh$ such that $L[\ga i]=L[\gao i]$, 
i=1,2. Then
\[
d(\ga 1 (0), \ga 2 (l_2)) \leq d(\gao 1 (0),\gao 2 (\lo 2)).
\]
\end{lemma}
\begin{thm}{\bf (Topogonov)}
Let $(M,g)$ be a  Hopf-Rinow manifold such that 
$H \leq K \leq \Delta$. Then
\begin{itemize}
\item [{\bf (A)}] let $\tr$ be a geodesic triangle in $M$.
Assume $\ga 1,$ $\ga 3$ are minimal geodesic and if $H >0$,
$l_2 \leq \frac{\pi}{\sqrt H}$. Then in $\mh,$ simply connected 
$2$-dimensional manifold of constant curvature $H,$ there exists
a geodesic triangle $\tro$ such that $l_i=\ov{l}_i$ and $\alfo 1 \leq \alf 1$,
$\alfo 2 \leq \alf 2$. Except in case $H>0$ and $l_2 =\frac{\pi}{\sqrt H},$
the triangle is uniquely determined.  
\item [{\bf (B)}] Let $\gga$ be a hinge in $M.$ Let 
$\ga 1$ be a minimal geodesic, and if $H>0$, $l_2 \leq \frac{\pi}{\sqrt H}$. 
Let $\ggao$ be a hinge in $\mh$ such that 
$l_i =\ov{l}_i,\ i=1,2$ e $\alpha=\ov{\alpha}$. Then
\[
d(\ga 1 (0), \ga 2 (0)) \leq d( \gao 1 (0), \gao 2 (0)) .
\]
\end{itemize}
\end{thm}
\dem $\ $ the proof consists of numbers of steps as in \cite{CE} 
page.43 . First of all, we recall briefly some facts of the proof in
\cite{CE}. Let $\gga$ be a hinge. We call this hinge small if
$\frac{1}{2} r=  max \ L[ \gamma_{i} ],\ i=1,2$ and $\e{ \ga 2 (0)}$ is an
embedding on $\bn{p}{r}$. Let  $\tr$ be a geodesic triangle. We call $\tr$ 
a {\em small triangle} if any hinge of $\tr$ is small. 
Let $\tr $ as in  {\bf (A)}. We say that
$\tr$ is thin if $( \gamma_1,\gamma_2,\alpha_3 )$ and 
$( \gamma_3,\gamma_2,\alpha_1 )$ are thin hinges, i.e. thin right hinge or 
thin obtuse hinge or thin acute hinge. We briefly describe the above
terminology.  

A thin right hinge is a hinge 
$\trh$ if the hypothesis of corollary \ref{cb} hold. 

Let $\gga$ be a hinge
with $\alpha > \frac{\pi}{2}$. Let $\trho$ be the corresponding hinge in
$\mhe$ with $L[\gao i] =L [\ga i],\ i=1,2.$ Let $\gao 3$ be a minimal 
geodesic from $\gao 2 (l_2)$ to $\gao 1 (0).$ 
Let $\ov{\sigma}:[0,l] \lra \mhe$ be a geodesic starting 
from $\gao 2 (0)$ such that
\[
\langle \fis{\ov{\sigma}}(0),\fis{\gao 2}(0) \rangle =0,\ 
\fis{\ov{\sigma}}(0)=\lambda_1 \fis{\gao 1}(l_1) \ + \ 
\lambda_2 \fis{\gao 2}(0),\ \lambda_i>0
\]
and let $\ov{\sigma}(l)$ be the first point of $\ov{\sigma}$ which lies on
$\gao 3.$  Let $\sigma$ be a geodesic in $M$
starting from $\ga 2 (0)$
with the same properties of $\ov{\sigma}$. We call  $\gga$ is 
a {\em thin obtuse hinge} if 
$(\ga 1,\sigma,\alpha-\frac{\pi}{2})$ is a small hinge and
$(\sigma,\ga 2, \frac{\pi}{2})$ is a thin right hinge.   

Let $\gga$ be a hinge with $\alpha< \frac{\pi}{2}$ and let $\ga 2 (l)$
be a point closest to $\ga 1 (0)$ .
Let 
\begin{center}
$\tau=\ga 2:[0,l] \lra M,\  $ \\
$\theta=\ga 2: [l,l_2]\lra M$
\end{center}
and let  $\sigma:[0,k] \longrightarrow M$ be a 
minimal geodesic from $\ga 1 (0)$ to
$\ga 2 (l_2)$. We call $\gga$ a {\em thin acute hinge} if 
$(\ga 1,\tau,\sigma)$ is a small triangle and  
$(\sigma,\theta,\frac{\pi}{2})$ is a small right hinge.

From step $(1)$ to  step $(7),$ in \cite{CE}, they essentially prove that 
{\bf (B)} holds for thin right hinges, thin obtuse hinge and thin acute
hinge. The same proofs work in our context since in any steps 
they use  Rauch and Berger theorems, and  the main corollaries, 
and the existence of at least a minimal geodesic joining any two points.

Now, we will prove it in general.   Given an arbitrary
hinge $\gga$ as in {\bf (B)}, fix $N$ and let
\[
\tau_{k,l}= \ga 2:[\frac{kl_2}{N},\frac{(k+l)l_2}{N}] \lra M,
\]
where $k,l$ are integers with $0 \leq k,l \leq N.$ Let $\sigma_k$ be 
the minimal geodesic from $\ga 1 (0)$ to $\ga 2 (\frac{k l_2}{N}).$ As
in \cite{CE} page. 48 we shall prove that any triangle 
$T_{k,l}=(\sigma_k,\tau_{k,l},\sigma_{k+l})$ is a geodesic triangle. If we 
prove that there exists $N$ such that any $T_{k,1}$ is thin we may continue
the proof as in \cite{CE}, proving our aim. 
If both $ H \leq 0$ and $\Delta$ are 
non positive the result follows easily while  
if $H \leq 0$ and $\Delta>0$ then it is 
enough to choose $N$ such that $(i)$ in lemma \ref{toro} holds. 
Then we shall assume $H>0$ and by Berger-Topogonov theorem we may suppose
also that $ d (\ga 1 (0), \ga 2 (t) ) < \frac{\pi}{\sqrt{H}} .$ Indeed, 
if $d (\ga (0), \ga 2 (t)= \frac{\pi}{\sqrt{H}}$ for some $t$ then the
manifold must be isometric to a sphere concluding our proof.
Using the compactness of $\ga 2$ there exists $\epsilon_o>0$ such that for 
every
$\epsilon \leq \epsilon_o$ there exists $\eta (\epsilon)$ such that 
for every $
s,t \in [0, l_2]$ we get
\[
d(\gamma_1(0),\gamma_2 (t)) \ + \ d(\gamma_1(0), \gamma_2(s)) \leq
\frac{2 \pi}{\sqrt{H - \epsilon}}\ - 5\eta(\epsilon).
\] 
Hence, for every $\epsilon \leq \epsilon_o,$ we choose  $N$ such 
that $L[ \ti{k,1}] \leq \eta(\epsilon)$ and once apply lemma \ref{toro}, 
comparing $M$ with $\mhe,$ on
$(\si k,\ti{k,1},\alf k)$  and  $(\si{k+l},\ti{k,i},\bi k ).$ Moreover, using
again the compactness of $\ga 2$ then there exists $r_o$ such that 
\[
\e {\ga 2  (t)}: \bn{\ga 2 (t)}{2r_o} \longrightarrow \bd{B}{\ga 2 (t)}{2r_o}
\]
is a diffeomorphism onto. Now it is easy to see that any geodesic triangle
$T_{k,l}$ is thin. QED
\begin{cor} \label{pippo}
Let $\tr$ be a geodesic triangle in a Hopf-Rinow manifold 
such that $0<H \leq K \leq \Delta$. Then the perimeter of
$\tr$ is at most $\frac{2 \pi }{\sqrt{H}}.$
\end{cor}
\begin{cor} \label{bt2}
Let $(M,g)$ be a Hopf-Rinow manifold with
$0<H \leq K \leq \Delta$. Assume that
$d(M)= \frac{\pi}{\sqrt H}$ and that there exists a point 
$p$ such that the image of the function $q \vai d(p,q)$ 
has $[0,\frac{\pi}{\sqrt H})$ as subset. Then
$M$ is isometric to
$S_{\frac{1}{H}}(\tpm \times \R).$
\end{cor}
\n
\dem $\ $ by Ekeland theorem (\cite{Ek}) there exists a sequence
$q_n,$ in $M$ such that
\[
d(p,q_n) \vai \frac{\pi}{\sqrt H}.
\]
and there exists a unique geodesic $\gamma_n,$ that we shall assume 
parameterized in $[0,1],$ from $p$ to $q_n$. Take the hinges
$(\ga n , \ga m, \alf{n,m}).$ Using Topogonov theorem there exists 
a hinge $(\gao n, \gao m, \alf{n,m})$ in $\mh$ such that
\[
d(q_n,q_m) \leq d(\gao n (1),\gao m (1)),
\]
Now,  $\gao n (1)$ converges then $q_n$ is a Cauchy sequence in $M.$ 
In particular there exists the limit $q$ of the sequence $q_n$ that it
satisfies $d(p,q)=d(M).$ Using Berger-Topogonov theorem,
$M$ is isometric to the sphere
$ S_{\frac{1}{\sqrt{H}}}(\tpm \times \R)$. QED 
\begin{thm} {\bf Sphere theorem}
Let $(M,g)$ be a Hopf-Rinow manifold such that
$\frac{4}{9}<\delta \leq K \leq 1$. Assume that
$i(M)\geq \pi$. Then $M$ is contractible and if $\H $ is 
separable then $M$ is diffeomorphic to $S(l_2)$.
\end{thm}
\n
\dem $\ $ since $\frac{4}{9} < \delta$   there exists 
$\epsilon >0$ such that
\[
\frac{\pi}{\sqrt{\delta}}=\frac{3}{2}(\pi \ - \ \epsilon).
\]
Using the fact $i(M)\geq \pi,$ there exist two points 
$p,q \in M$ such that $d(p,q)=\pi \ - \ \epsilon$. We claim that
$M$ is covered by the following geodesic balls
\[
M= \ov{\bd{\cal B}{p}{\pi-{\epsilon}}} \cup
\ov{\bd{\cal B}{q}{\pi-{\epsilon}}}.
\]
Let $r\in M$ such that $d(p,r) \geq \pi \ - \ \epsilon$. Using  corollary  
\ref{pippo}, we have
\[
d(q,r)\leq 2(\frac{3}{2}(\pi \ - \ \epsilon)) \ - \ 2(\pi \ - \ \epsilon)
      =\pi \ - \ \epsilon.
\]
Now, we shall conclude our proof as in the Sphere Rauch theorem.
QED

\n
Leonardo Biliotti \\
Universit\`a degli studi di Firenze \\
Dipartimento di Matematica e delle Applicazioni all'Architettura \\
Piazza Ghiberti 27 - Via dell'Agnolo 2r - 50132 Firenze (Italy) \\
e-mail:{\tt biliotti@math.unifi.it} \\

\begin{thebibliography}{10}
\bibitem{AMR}
{\sc {A}braham, R., Marsden, J.~E., and Ratiu, T.}
\newblock {\em Manifolds, tensor analysis, and application}.
\newblock Applied Mathematical Sciences 75, Springer-Verlang, New York, 1988.
%
\bibitem{An}
{\sc Anderson, L.}
\newblock {\em The {B}onnet-{M}yer theorem is true for {R}iemannian {H}ilbert
  manifold}.
\newblock Math. Scand. 58 (1986), 236--238.
%
\bibitem{At1}
{\sc Atkin, C.~J.}
\newblock {\em Geodesic and metric completeness in infinite dimensions}.
\newblock Hokkaido Mathematical Journal  vol 26 (1997), 1--61.
%
\bibitem{At}
{\sc Atkin, C.~J.}
\newblock {\em The {H}opf-{R}inow theorem is false in infinite dimension}.
\newblock  Bull. London Math. Soc. 7 (1975), 261--266.
%
\bibitem{BE}
{\sc Bessega, C.}
\newblock {\em Every infinite-dimensional {H}ilbert space is 
diffeomorphic with its unit sphere}.
\newblock Bull. Acad. Polon. Sci XIV (1966), 27--31.
%
\bibitem{B}
{\sc Biliotti, L.}
\newblock{\em Properly discontinuous isometric actions on the unit
sphere of infinite dimensional Hilbert spaces}
\newblock{preprint} 
%
\bibitem{BMT}
{\sc Biliotti, L., Mercuri, F., and Tausk, D.}
\newblock {\em Note on tensor fields in {H}ilbert spaces}.
\newblock{Academia Brasileira de Ci\^encies (2002) 74 (2),  207-210.}
%
\bibitem{KB}
{\sc {B}urghelea, D., and Kuiper, N.}
\newblock {\em {H}ilbert manifolds}.
\newblock  Annals of Mathematics 90  (1968), 379--417.
%
\bibitem{CE}
{\sc Cheeger, J., and Ebin, D.}
\newblock {\em Comparison theorem in {R}iemannian geometry}.
\newblock North--Holland, Amsterdam, 1975.
%
\bibitem{Ee}
{\sc Eells, J.}
\newblock {\em A setting for global analysis}.
\newblock Bull. Amer. Math. Soc 72  (1966), 751--807.
%
\bibitem{Ek}
{\sc Ekeland, I.}
\newblock {\em The {H}opf-{R}inow theorem in infinite dimension}.
\newblock Journal of Differential Geometry 13 (1978), 287--301.
%
\bibitem{El}
{\sc El\'{\i}ason, H.}
\newblock {\em Condiction (C) and geodesic on {S}obolev manifolds}.
\newblock Bull. Amer. Math. Soc. 77, n 6 (1971), 1002--1005.
%
\bibitem{Gr2}
{\sc Grossman, N.}
\newblock {\em {G}eodesic on {H}ilbert manifolds}.
\newblock PhD thesis, University of Minnesota, 1964.
%
\bibitem{Gr}
{\sc Grossman, N.}
\newblock {\em {H}ilbert manifolds without epiconjugates points}.
\newblock Proc. of Amer. Math. Soc 16 (1965), 1365--1371.
%
\bibitem{Ka}
{\sc Karcher, H.}
\newblock {\em Riemannian center of mass and mollifer smoothing}.
\newblock  Comm. on Pure and Appl. Math. XXX (1977), 509--541.
%
\bibitem{Kl}
{\sc Klingemberg, W.}
\newblock {\em Riemannian geometry}.
\newblock De Gruyter studies in Mathemathics, New York, 1982.
%
\bibitem{KN}
{\sc Kobayashi, S., and Nomizu, K.}
\newblock {\em Foundations of differential geometry}.
\newblock vol I, Interscience Wiley, New York, 1963.
%
\bibitem{La}
{\sc Lang, S.}
\newblock {\em Differential and {R}iemannian manifolds}.
\newblock Third Ediction. Graduate Text em Mathematics, 160, Springer-Verlang,
New York, 1996.
%
\bibitem{Mca}
{\sc McAlphin, J.}
\newblock {\em Infinite dimensional manifolds and {M}orse theory}.
\newblock PhD thesis, University of Columbia, 1965.
%
\bibitem{Mi}
{\sc Misiolek, G.}
\newblock {\em {T}he exponential map on the free loop space if {F}redholm}.
\newblock Geom.  Funct. Anal. 7 (1997), 1--17.
%
\bibitem{Pa}
{\sc Palais, R.}
\newblock {M}orse {t}heory {o}n {H}ilbert {m}anifolds.
\newblock {\em Topology 2\/} (1963), 299--340.
%
\bibitem{Pa2}
{\sc Palais, R.}
\newblock {\em {H}omotopy {t}heory {o}f {i}nfinite {d}imensional {m}anifolds}.
\newblock  Topology 5 (1966), 1--16.
, Inc., New York-Amsterdam, 1968.
%
\bibitem{Ra}
{\sc Rauch, H.~E.}
\newblock  {\em A contribution of differential geometry in the large}.
\newblock Annals of Mathematics 54 (1951), 38--55.
%
\bibitem{Sa}
{\sc Sakai, T.}
\newblock {\em Riemannian geometry}.
\newblock Trans. Math. Monog. vol 149. AMS, New York, 1996.
%
\end{thebibliography}
\end{document}